# Null eikonal helices and Null eikonal slant helices in the 4-dimensional Lorentzian manifold

**Evren Zıplar**
*Çankırı Karatekin University, Faculty of Science, Department of Mathematics, Çankırı, Turkey*
E-mail: evrenziplar@karatekin.edu.tr

**Abstract**

In this paper, we define the notion of eikonal helix and eikonal slant helix for null curves in the 4-dimensional Lorentzian manifold $M_1^4$ and give a characterization for the null curve to be the null eikonal helix. Moreover, we indicate an important relation between the null eikonal helix and null eikonal slant helix and find the axis of the null eikonal helix. We obtain some relationships between the curvatures of these curves.

## 1. Introduction

In geometrical structures of optical lenses, the eikonal function is an essential factor in the microwave frequency range. Practically, it is utilized as an option to Snell's law which is a formula used to express the relationship between the angles of incidence and refraction. Moreover, the geometrical equations of ray trajectories in media with certain shapes of refractive index variation is occurred by the eikonal equation generated by the eikonal function [1].

On the other hand, 4-dimensional Lorenzian manifolds as compared to the Euclidean manifolds or lower dimensional Lorenzian manifolds are much more exciting. In relativity theory, these manifolds have a significant act. In this theory, a unit speed future-pointing timelike curve in a space-time which is a geometrical model that unifies space and time is thought as a place of a material particle in the space-time. The unit-speed parameter of this curve is considered as the appropriate time of a material particle. Similar to the timelike curves, a future-pointing null geodesic is imagined as the place of a lightlike particle. More generally, the working of null curves has its own geometric concern. Because the other curves, such as spacelike and timelike, can be worked by a similar way to that examined in positive definite Riemannian geometry. Futhermore, null curves as compared to spacelike and timelike curves have different features. In consequence of the growing importance of null curves in geometrical physics, some particular subjects of the curve theory have been worked for null curves by many scientists dealing with differential geometry. Ferrandez, Gimenez and Lucas have introduced and examined null helices in Lorentzian space forms [4]. Honda and Inoguchi have given a characterization of null cubics [5]. Thereafter, Karadağ H.B. and Karadağ M. have investigated null generalized slant helices in Lorentzian space [6]. Surely, these topics are investigated in 3-dimensional Lorenzian space.

This paper is organized as follows. In section 2, we will give some basic notions and introduce the new kinds of null curves. In section 3, we give examples of such curves. In section 4, we will prove a Theorem which characterizes the null eikonal helix. In section 5, we obtain the axis of null $f$-eikonal helix and give a relation between the null $f$-eikonal helix and null $f$-eikonal slant helix. Finally, in section 6, we give a characterization of null $f$-eikonal slant helix.

## 2. Basic properties

In this section, we give the basic notions in Lorentzian manifolds. Also, we introduce two new Lorentzian curves which are called null eikonal helix and null eikonal slant helix in Definition 3.

**Definition 1.** ([3]) Let $g$ denote the metric in 4-dimensional Lorentzian manifold $M_1^4$, where $g(a,b) = -a_1b_1 + a_2b_2 + a_3b_3 + a_4b_4$ for the vectors $a = (a_1,a_2,a_3,a_4)$, $b = (b_1,b_2,b_3,b_4) \in TM_1^4$. A curve $\alpha(t)$ in $M_1^4$ is called a null curve if $g(\alpha'(t),\alpha'(t)) = 0$ and $\alpha'(t) \neq 0$ for all $t$. We note that a null curve $\alpha(t)$ in $M_1^4$ satisfies $g(\alpha''(t),\alpha''(t)) \geq 0$. We say that a null curve $\alpha(t)$ in $M_1^4$ is parametrized by pseudo-arc if $g(\alpha''(t),\alpha''(t)) = 1$.

Let us say that a null curve $\alpha(t)$ in $M_1^4$ with $g(\alpha''(t),\alpha''(t)) \neq 0$ is a Cartan curve if $\{\alpha'(t),\alpha''(t),\alpha'''(t)\}$ is linearly independent for any $t$. For a Cartan curve in $M_1^4$ with pseudo-arc parameter $t$, there exists a unique Frenet frame $\{\xi = \alpha', N, W_1, W_2\}$ such that

$$\dot{\xi} = \nabla_\xi \xi = W_1,$$
$$\dot{N} = \nabla_\xi N = \sigma_1 W_1 + \sigma_2 W_2,$$
$$\dot{W_1} = \nabla_\xi W_1 = -\sigma_1 \xi - N,$$
$$\dot{W_2} = \nabla_\xi W_2 = -\sigma_2 \xi,$$

where $N$ is null, $g(\xi, N) = 1$, $\{\xi, N\}$ and $\{W_1, W_2\}$ are orthogonal, $\{W_1, W_2\}$ is orthonormal, $\nabla$ is Levi-Civita connection in $M_1^4$. The functions $\{\sigma_1, \sigma_2\}$ are called the Cartan curvatures of $\alpha$ [2,4].

**Definition 2.** Let $M_1^4$ be a Lorentzian manifold with metric $g$ and let $\alpha$ be a curve in $M_1^4$. For the function $f : M_1^4 \to \mathbb{R}$, it is said that $f$ is eikonal along $\alpha$ if $\|\nabla f\|$ is constant along $\alpha$, where $\nabla f$ is gradient of $f$, i.e. $df(X) = g(\nabla f, X) = X(f)$ for all $X \in \mathfrak{X}(M_1^4)$.

**Lemma 1.** ([7]) Let $(M,g)$ be a pseudo-Riemannian manifold and $\nabla$ be the Levi-Civita connection of $M$. The Hessian $H^f$ of a $f \in F(M)$ is the symmetric tensor field such that

$$H^f(X,Y) = g(\nabla_X (\text{grad} f), Y),$$

where $F(M)$ shows the set of differentiable functions defined on $M$.

From Lemma 1, we have the following corollary.

**Corollary 1.** ([8]) The Hessian $H^f$ of a $f \in F(M)$ is zero, i.e., $H^f = 0$ if and only if $\nabla f$ is parallel in $M$.

**Definition 3.** Let $M_1^4$ be a Lorentzian manifold with metric $g$ and $\alpha(t)$ be a null curve with the Frenet frame $\{\xi = \alpha', N, W_1, W_2\}$. Let $f : M_1^4 \to \mathbb{R}$ be an eikonal function along the curve $\alpha$, i.e. $\|\nabla f\|$ is constant along $\alpha$. Then, we define the followings:

i) If the function $g(\nabla f, \xi)$ is a non-zero constant along $\alpha$, then $\alpha$ is called a null $f$-eikonal helix. And, $\nabla f$ is called the axis of the null $f$-eikonal helix $\alpha$.

ii) If the function $g(\nabla f, N)$ is a non-zero constant along $\alpha$, then $\alpha$ is called a null $f$-eikonal slant helix. And, $\nabla f$ is called the axis of the null $f$-eikonal slant helix $\alpha$.

## 3. Examples in Lorentzian manifold $\mathbb{R}_1^4$

In this section, we give examples of null eikonal helix and null eikonal slant helix.

**Example 1.** Consider the null curve $\alpha(t)$ in $\mathbb{R}_1^4$ given by the parametrization

$$\alpha(t) = -\frac{1}{\sqrt{2}}(\sinh t, \cosh t, \cos t, \sin t)$$

with pseudo-arc parameter $t$. Then the tangent vector $\xi$ is given by

$$\xi = -\frac{1}{\sqrt{2}}(\cosh t, \sinh t, -\sin t, \cos t).$$

On the other hand, let consider the function

$$f : \mathbb{R}_1^4 \to \mathbb{R}, \quad (x_1, x_2, x_3, x_4) \to f(x_1, x_2, x_3, x_4) = x_1 x_2$$

If we compute $\nabla f$, we find that $\nabla f = (x_2, x_1, 0, 0)$. Then, we have

$$\|\nabla f\| = \sqrt{|x_1^2 - x_2^2|} = \frac{1}{\sqrt{2}} = const.$$

along curve $\alpha$. That is, $f$ is an eikonal function along $\alpha$. So, we can easily see that $g(\nabla f, \xi) = -\frac{1}{2} = const.$, where $\nabla f = -\frac{1}{\sqrt{2}}(\cosh t, \sinh t, 0, 0)$ along $\alpha$. Finally, $\alpha$ is a null $f$-eikonal helix with the axis $\nabla f$.

**Example 2.** Consider the null curve $\alpha(t)$ in $\mathbb{R}_1^4$ given by the parametrization

$$\alpha(t) = -\left(\frac{1}{6}t^3 + t,\ \frac{1}{2}t^2,\ t,\ \frac{1}{6}t^3\right)$$

with pseudo-arc parameter $t$. Moreover, we can choose the Frenet frame elements of the curve $\alpha$ as followings:

$$N = (1,0,0,1)$$
$$W_1 = (-t,-1,0,-t)$$
$$W_2 = (-1,0,-1,-1),$$

where the tangent vector is $\xi = -\left(\frac{t^2}{2}+1,\ t,\ 1,\ \frac{t^2}{2}\right)$. In fact, we can easily see that $N$ is null, $g(\xi, N) = 1$, $\{\xi, N\}$ and $\{W_1, W_2\}$ are orthogonal, $\{W_1, W_2\}$ is orthonormal. Moreover, they verify the Frenet frame. Let consider the projection map

$$f : \mathbb{R}_1^4 \to \mathbb{R},\quad (x_1, x_2, x_3, x_4) \to f(x_1, x_2, x_3, x_4) = x_4.$$

If we compute $\nabla f$, we obtain that $\nabla f = (0,0,0,1)$. Then, we have $\|\nabla f\| = 1 = \text{constant}$ along curve $\alpha$. Then, $f$ is an eikonal function along $\alpha$. Furthermore, we can see that $g(\nabla f, N) = 1$. That is, $\alpha$ is a null $f$-eikonal slant helix with the axis $\nabla f$.

### 4. Characterization theorem

In this section, we give a characterization for a null curve to be a null $f$-eikonal helix.

**Theorem 1.** Let $\alpha(t)$ be a null curve with the tangent vector field $\xi = \alpha'$ in $M_1^4$. Let $f : M_1^4 \to \mathbb{R}$ be an eikonal function along the curve $\alpha$. Then, $\alpha$ is a null $f$-eikonal helix if and only if $f$ is a linear function along $\alpha$.

**Proof.** Let $\alpha$ be a null $f$-eikonal helix. Due to the fact that $f$ is an eikonal function along $\alpha$, $\|\nabla f\|$ is constant along $\alpha$. On the other hand, it is known that $X(f) = g(\nabla f, X)$ for all $X \in \mathfrak{X}(M_1^4)$. So, we have

$$g(\nabla f, \xi) = \xi(f) = \frac{d}{dt}(f \circ \alpha)$$

along $\alpha$. Since $\alpha$ is a null $f$-eikonal helix, $g(\nabla f, \xi)$ is constant along $\alpha$. Thus, we obtain

$$\frac{d}{dt}(f \circ \alpha) = const.$$

along $\alpha$. In other words, $f$ is a linear function along $\alpha$.

Conversely, we suppose that $f$ is a linear function along $\alpha$. Clearly,

$$\frac{d}{dt}(f \circ \alpha) = const.$$

Therefore, we have $g(\nabla f, \xi)$ is constant along $\alpha$. This completes the proof.

### 5. The axis of the null eikonal helix

Here, we find the axis of a null $f$-eikonal helix and give an important characterization of null $f$-eikonal helix. Moreover, we point out a relation between the null $f$-eikonal helix and null $f$-eikonal slant helix.

**Theorem 2.** Let $\alpha(t)$ be a null curve in $M_1^4$ with non-zero Cartan curvatures $\sigma_1, \sigma_2$ and let assume that $f: M_1^4 \to \mathbb{R}$ be an eikonal function along $\alpha$ and $H^f = 0$. If $\alpha(t)$ is a null $f$-eikonal helix in $M_1^4$, then the followings hold:

i) $\int \sigma_2 dt - \dfrac{\sigma_1'}{\sigma_2}$ is a constant function.

ii) The axis of the null eikonal helix is $\nabla f = c\left[(-\sigma_1)\xi + N + \left(-\int \sigma_2 dt\right)W_2\right]$, where $g(\nabla f, \xi) = c$ is a non-zero constant.

iii) $\sigma_1$ is non-constant.

**Proof. i)** Let $\alpha(t)$ be a null $f$-eikonal helix curve in $M_1^4$ with axis $\nabla f$. Then, there exist smooth functions $a_1 = a_1(t)$, $a_2 = a_2(t)$, $a_3 = a_3(t)$ and $a_4 = a_4(t)$ of pseudo-arc parameter $t$ such that

$$\nabla f = a_1 \xi + a_2 N + a_3 W_1 + a_4 W_2, \tag{1}$$

where $\{\xi, N, W_1, W_2\}$ is a basis of $\mathfrak{X}(\mathbb{R}_1^4)$. From (1), we have

$$\begin{cases} g(\nabla f, \xi) = a_2 = c = const., \\ g(\nabla f, N) = a_1, \\ g(\nabla f, W_1) = a_3, \\ g(\nabla f, W_2) = a_4. \end{cases} \tag{2}$$

By taking derivative in each part of equations in (2) in the direction $\xi$ in $M_1^4$ and using Frenet frame, we get

$$\begin{cases} a_3 = 0, \\ \sigma_2 a_4 = a_1', \\ c\sigma_1 = -a_1, \\ c\sigma_2 = -a_4'. \end{cases} \tag{3}$$

Using the second and third equations of the system (3), we have

$$a_4' = -c\left(\frac{\sigma_1'}{\sigma_2}\right)'. \tag{4}$$

From (4) and the last equation of the system (3), it follows

$$\sigma_2 - \left(\frac{\sigma_1'}{\sigma_2}\right)' = 0.$$

In other words, $\int \sigma_2 dt - \frac{\sigma_1'}{\sigma_2}$ is a constant function.

**ii)** From the system (3) and (1), the axis of the null eikonal helix curve is obtained as $\nabla f = c\left[(-\sigma_1)\xi + N + \left(-\int \sigma_2 dt\right)W_2\right]$, where $g(\nabla f, \xi) = c$ is a non-zero constant.

**iii)** Let $\sigma_1$ be a constant function. Then, from the system (3), we get $\sigma_2 = 0$. But, we have a contradiction since all curvatures are nowhere zero. This completes the proof.

The above theorem has the following corollary.

**Corollary 2.** Let $\alpha(t)$ be a null curve in $M_1^4$ with non-zero Cartan curvatures $\sigma_1, \sigma_2$ and let assume that $f: M_1^4 \to \mathbb{R}$ be an eikonal function along $\alpha$ and $H^f = 0$. If $\alpha(t)$ is a null $f$-eikonal helix in $M_1^4$, then $\alpha(t)$ can not be a null $f$-eikonal slant helix in $M_1^4$.

**Proof.** Let $\alpha(t)$ be a null $f$-eikonal helix in $M_1^4$. Then, from Theorem 2, the axis of the helix $\alpha(t)$ is

$$\nabla f = c\left[(-\sigma_1)\xi + N + \left(-\int \sigma_2 dt\right)W_2\right], \tag{5}$$

where $g(\nabla f, \xi) = c$ is a non-zero constant. By using the equation (5), we have

$$g(\nabla f, N) = -c\sigma_1. \tag{6}$$

From Theorem 2, $\sigma_1$ is non-constant, therefore $g(\nabla f, N)$ can not be a constant. This completes the proof.

**Theorem 3.** Let $\alpha(t)$ be a null curve in $M_1^4$ with non-zero Cartan curvatures $\sigma_1, \sigma_2$ and let assume that $f: M_1^4 \to \mathbb{R}$ be an eikonal function along $\alpha$ and $H^f = 0$. If $\alpha(t)$ is a null $f$-eikonal helix in $M_1^4$, then $\det\left(\nabla_\xi \xi,\ \nabla_\xi^2 \xi,\ \nabla_\xi^3 \xi,\ \nabla_\xi^4 \xi\right) = 0$ holds.

**Proof.** Let $\alpha(t)$ be a null $f$-eikonal helix curve in $M_1^4$. Hence, from Frenet frame, we have the followings

$$\begin{cases} \nabla_\xi \xi = W_1, \\ \nabla^2_\xi \xi = -\sigma_1 \xi - N, \\ \nabla^3_\xi \xi = \left(-\sigma_1'\right)\xi + (-2\sigma_1)W_1 + (-\sigma_2)W_2, \\ \nabla^4_\xi \xi = \left(-\sigma_1'' + \sigma_2^2 + 2\sigma_1^2\right)\xi + (2\sigma_1)N + \left(-3\sigma_1'\right)W_1 + \left(-\sigma_2'\right)W_2. \end{cases} \qquad (7)$$

From (7), we get $\det\left(\nabla_\xi \xi, \nabla^2_\xi \xi, \nabla^3_\xi \xi, \nabla^4_\xi \xi\right) = \sigma_2^3 - \left(\sigma_1'' \sigma_2 - \sigma_2' \sigma_1'\right)$. By using Theorem 2 (i), it can be seen that $\left(\sigma_1'' \sigma_2 - \sigma_2' \sigma_1'\right) = \sigma_2^3$. Finally, $\det\left(\nabla_\xi \xi, \nabla^2_\xi \xi, \nabla^3_\xi \xi, \nabla^4_\xi \xi\right) = 0$. This completes the proof.

## 6. Characterization for the null eikonal slant helix

Here, we give a characterization of the null eikonal slant helix in Theorem 4. Then, we investigate a particular case of the Theorem.

**Theorem 4.** Let $\alpha(t)$ be a null curve in $M_1^4$ with Cartan curvatures $\sigma_1, \sigma_2$ and let assume that $f : M_1^4 \to \mathbb{R}$ be an eikonal function along $\alpha$ and $H^f = 0$. If $\alpha(t)$ is a null $f$-eikonal slant helix in $M_1^4$, then the following system is satisfied:

$$\begin{cases} a_2 \sigma_1 + a_3 \sigma_2 = 0, \\ a_1' - a_2 = 0, \\ a_2' + a_1 \sigma_1 + c = 0, \\ a_3' + a_1 \sigma_2 = 0, \end{cases} \qquad (8)$$

where $a_1 = g(\nabla f, \xi)$, $a_2 = g(\nabla f, W_1)$, $a_3 = g(\nabla f, W_2)$ and $g(\nabla f, N) = c$ is a non-zero constant.

**Proof.** Let $\alpha(t)$ be a null $f$-eikonal slant helix in $M_1^4$ with axis $\nabla f$. Then, there exist smooth functions $a_1 = a_1(t)$, $a_2 = a_2(t)$, $a_3 = a_3(t)$ and $a_4 = a_4(t)$ of pseudo-arc parameter $t$ such that

$$\nabla f = c\xi + a_1 N + a_2 W_1 + a_3 W_2 , \qquad (9)$$

where $g(\nabla f, N) = c$ is non-zero constant. From Corollary 1, $\nabla f$ is parallel in $M_1^4$, i.e., $\nabla_\xi (\nabla f) = 0$ along $\alpha$. Hence, taking derivative in each part of (9) in the direction $\xi$ and using the Frenet frame, we obtain

$$(-a_2 \sigma_1 - a_3 \sigma_2)\xi + \left(a_1' - a_2\right)N + \left(a_2' + a_1 \sigma_1 + c\right)W_1 + \left(a_3' + a_1 \sigma_2\right)W_2 = 0, \qquad (10)$$

which gives the following system

$$\begin{cases} a_2\sigma_1 + a_3\sigma_2 = 0, \\ a_1' - a_2 = 0, \\ a_2' + a_1\sigma_1 + c = 0, \\ a_3' + a_1\sigma_2 = 0. \end{cases}$$

This completes the proof.

The above Theorem has the following corollary.

**Corollary 3.** In Theorem 4, if $\sigma_1 = 0$, then the axis of the null $f$-eikonal slant helix is given by

$$\nabla f = c\xi + \left[\left(-\frac{c}{2}\right)t^2 + mt + n\right]N + (-ct + m)W_1 + kW_3, \tag{11}$$

where $c \neq 0, m, n, k$ are constants.

**Proof.** Let $\sigma_1 = 0$. Then, by using the second and the third equations of system (8), we obtain the second order differential equation

$$a_1''(t) = -c. \tag{12}$$

The solution of (12) is

$$a_1(t) = \left(-\frac{c}{2}\right)t^2 + mt + n, \tag{13}$$

where $c, m$ and $n$ are constants. By using the second equation of (8) and taking the derivative of both sides of (13), we have

$$a_2(t) = -ct + m. \tag{14}$$

On the other hand, from the first and the fourth equations of (8), $\sigma_2$ must be zero since $a_1$ is different from zero. Thus, by using the fourth equation of (8), we get

$$a_3 = k = const. \tag{15}$$

Finally, by using (13), (14) and (15), the axis of the $f$-eikonal slant helix is obtained as follows

$$\nabla f = c\xi + \left[\left(-\frac{c}{2}\right)t^2 + mt + n\right]N + (-ct + m)W_1 + kW_3,$$

where $c \neq 0, m, n, k$ are constants and $\sigma_1 = \sigma_2 = 0$. This completes the proof.

**Example 3.** We consider the curve $\alpha(t)$ given in example 2. Then, from Frenet frame, we can see that $\sigma_1 = \sigma_2 = 0$. In the equality (11), if we take $c=1$, $k=-1$, $m=n=0$, we find the axis of $\alpha(t)$ as

$$\nabla f = \xi + \left(-\frac{t^2}{2}\right)N + (-t)W_1 - W_2$$
$$\nabla f = (0,0,0,1).$$